\documentclass[11pt]{amsart}

\usepackage{amscd}
\usepackage{amsmath}
\usepackage{graphicx}
\usepackage{amsfonts}
\usepackage{amssymb}
\begin{document}
\title[Strong skew commutativity preserving  maps]{
Strong skew commutativity preserving  maps on von Neumann algebras}

\author{Xiaofei Qi}
\address[Xiaofei Qi]{
Department of Mathematics, Shanxi University, Taiyuan 030006, P. R.
China.} \email{qixf1980@126.com}

\author{Jinchuan Hou}
\address[Jinchuan Hou]{Department of
Mathematics, Taiyuan University of Technology, Taiyuan 030024, P. R.
of China} \email{jinchuanhou@yahoo.com.cn}

\thanks{{ This work is partially supported by the National Natural Science Foundation of
China (11171249, 11101250).}}

\begin{abstract}

Let ${\mathcal M}$ be a von Neumann algebra without central summands
of type $I_1$.  Assume that $\Phi:{\mathcal M}\rightarrow {\mathcal
M}$ is a  surjective map. It is shown that $\Phi$ is strong skew
commutativity preserving (that is, satisfies
$\Phi(A)\Phi(B)-\Phi(B)\Phi(A)^*=AB-BA^*$ for all $A,B\in{\mathcal
M}$) if and only if there exists some self-adjoint element $Z$ in
the center of ${\mathcal M}$ with $Z^2=I$ such that $\Phi(A)=ZA$ for
all $A\in{\mathcal M}$. The strong skew commutativity preserving
maps on prime involution rings and prime involution algebras are
also characterized.

\end{abstract}
\maketitle

\section{Introduction}

Let $\mathcal R$ be a *-ring. For $A,B\in{\mathcal R}$, denote by
$[A,B]_*=AB-BA^*$ the skew Lie product. This product $AB-BA^*$ is
found playing an  important role in some research topics.  Let $B$
be a fixed element in $\mathcal R$. The additive map $\Phi:{\mathcal
R}\rightarrow{\mathcal R}$ defined by $\Phi(A)=AB-BA^*$ for all
$A\in{\mathcal R}$ is a Jordan
*-derivation, that is, it satisfies $\Phi(A^2)=\Phi(A)A^*+A\Phi(A)$. The
notion of Jordan *-derivations arose naturally in $\check{S}$emrl's
work \cite{S1,S2} investigating the problem of representing
quadratic functionals with sesquilinear functionals.  Motivated by
the theory of rings (and algebras) equipped with a Lie product
$[T,S]=TS-ST$ or a Jordan product $T\circ S=TS+ST$, Moln\'{a}r in
\cite{M} studied the skew Lie product   and  gave a characterization
of ideals in ${\mathcal B}(H)$ in terms of the skew Lie product. It
is shown in \cite{M} that, if ${\mathcal N}\subseteq {\mathcal
B}(H)$ is an ideal, then ${\mathcal N}={\rm span}\{AB-BA^*:
A\in{\mathcal N},B\in{\mathcal B}(H)\}={\rm span}\{AB-BA^*:
A\in{\mathcal B}(H),B\in{\mathcal N}\}$. In particular, every
operator in ${\mathcal B}(H)$ is a finite sum of $AB-BA^*$ type
operators. Later, Bre$\check{s}$ar and Fon$\check{s}$ner \cite{B}
generalized the above results in \cite{M} to rings with involution
in different ways.

\if false if $H$ is a real or complex Hilbert space with $\dim H>1$,
then a subspace $\mathcal N$ of ${\mathcal B}(H)$ is an ideal if and
only if $AB-BA^*\in{\mathcal N}$ for $A\in{\mathcal B}(H)$ and
$B\in{\mathcal N}$;  if $\dim H$ is an odd natural number, then
${\mathcal N}={\mathcal B}(H)$; and\fi

Recall that a  map $\Phi:{\mathcal R}\rightarrow{\mathcal R}$ is
skew commutativity preserving  if, for any $A$, $B\in {\mathcal R}$,
$[A,B]_*=0$ implies $[\Phi(A),\Phi(B)]_*=0$. The problem of
characterizing linear (or additive) bijective maps preserving skew
commutativity had been studied intensively on various  algebras (see
\cite{CFL, CH} and the references therein). More specially, we say
that a map $\Phi:{\mathcal R}\rightarrow{\mathcal R}$ is strong skew
commutativity preserving (briefly, SSCP) if
$[\Phi(A),\Phi(B)]_*=[A,B]_*$ for all $A,B\in{\mathcal R}$. SSCP
maps are also called strong skew Lie product preserving maps in
\cite{CP}. We prefer to the terminology SSCP because many authors
call the maps satisfying $[\Phi(A),\Phi(B)]=[A,B]$ strong
commutativity preserving maps. It is obvious that a strong skew
commutativity preserving map must be skew commutativity preserving,
but the inverse is not true generally. In \cite{CP}, Cui and Park
proved that, if $\mathcal R$ is a factor von Neumann algebra, then
every SSCP surjective map $\Phi$ on ${\mathcal R}$ has the form
$\Phi(A)=\Psi(A)+h(A)I$ for every $A\in{\mathcal R}$, where
$\Psi:{\mathcal R}\rightarrow{\mathcal R}$ is a linear bijective map
satisfying $\Psi(A)\Psi(B)-\Psi(B)\Psi(A)^*=AB-BA^*$ for
$A,B\in{\mathcal R}$ and $h$ is a real functional on $\mathcal R$
with $h(0)=0$; particularly, if $\mathcal R$ is a type $I$ factor,
then $\Phi(A)=cA+h(A)I$ for every $A\in{\mathcal R}$, where
$c\in\{-1,1\}$.

In the present paper, we show further that, in the above result,
$h=0$ and $\Psi (A)=A$ for all $A$ or $\Psi(A)=-A$ for all $A$. In
fact, the purpose of the present paper is to give a characterization
of the SSCP maps on prime  *-rings or on general von Neumann
algebras without central summands of type $I_1$. And the improvement
of the above result   in \cite{CP} is an immediate consequence of
our results.

Before embarking on the main results, we need some notations. Let
$\mathcal A$ be a *-ring. Denote by  ${\mathcal Z}({\mathcal A})$
the center of $\mathcal A$. Observe that $Z\in{\mathcal Z}({\mathcal
A})$ if and only if $Z^*\in{\mathcal Z}({\mathcal A})$. ${\mathcal
Z}_S({\mathcal A})=\{A\in{\mathcal Z}({\mathcal A}): A=A^*\}$. An
element $A\in{\mathcal A}$ is symmetric (respectively, skew
symmetric) if $A^*=A$ (respectively, if $A^*=-A$).
  If
$\mathcal A$ is a unital
*-algebra with unit $I$ over a field $\mathbb F$, let $\mathcal S$ (resp. $\mathcal
K$) be the set of its symmetric (resp. skew symmetric) elements.
Then every element $A\in{\mathcal A}$ can be written as $A=S+K$,
where $S\in{\mathcal S}$ and $K\in{\mathcal K}$. Moreover, this
decomposition is unique. We say that the involution $*$ is of the
first kind if ${\mathbb F}I\subseteq{\mathcal S}$; equivalently, $*$
is an $\mathbb F$-linear map. Otherwise, $*$ is said to be of the
second kind (see \cite{Br}). Recall that a ring $\mathcal A$ is
prime if, for any $A,\ B\in{\mathcal A}$, $A{\mathcal A}B=0$ implies
$A=0$ or $B=0$.  A von Neumann algebra $\mathcal M$ is a subalgebra
of some ${\mathcal B}(H)$, the algebra of all bounded linear
operators acting on a complex Hilbert space $H$, which satisfies the
double commutant property: ${\mathcal M}^{\prime\prime}=\mathcal M$,
where ${\mathcal M}^\prime =\{T: T\in{\mathcal B}(H)\ \mbox{and}\
TA=AT\ \forall A\in{\mathcal M}\}$ and ${\mathcal
M}^{\prime\prime}=\{{\mathcal M}^\prime\}^\prime$. For the theory of
von Neumann algebras, ref. \cite{KR}.

The rest part of this paper is organized as follows. In Section 2,
we discuss the question in a pure algebraic setting. Let $\mathcal
A$ be a unital prime
*-ring. Assume that $\Phi: {\mathcal A}\rightarrow {\mathcal A}$ is
a general strong skew commutativity preserving surjective map. We
show that, if $\mathcal A$ contains a nontrivial symmetric
idempotent element, then $\Phi$ has the form $\Phi(A)=A+f(A)$ for
all $A\in{\mathcal A}$, or $\Phi(A)=-A+f(A)$ for all $A\in{\mathcal
A}$, where $f:{\mathcal A}\rightarrow {\mathcal Z}_S({\mathcal A})$
is an arbitrary map (Theorem 2.1). Particularly, if $\mathcal A$ is
a
*-algebra and if the involution $*$ is of the second kind, then
$\Phi:{\mathcal A}\rightarrow{\mathcal A}$ is strong skew
commutativity preserving if and only if $\Phi(A)=A$ for all
$A\in{\mathcal A}$, or $\Phi(A)=-A$ for all $A\in{\mathcal A}$
(Theorem 2.10). As an application, we obtain a characterization of
strong skew commutativity preserving general surjective maps on
factor von Neumann algebras (Theorem 2.11) which improves the main
results in \cite{CP}. Section 3 is devoted to characterizing  the
strong skew commutativity maps on  general von Neumann algebras. We
prove that, if $\mathcal A$ is a von Neumann algebra without central
summands of type $I_1$, then a map $\Phi:{\mathcal
A}\rightarrow{\mathcal A}$ is strong skew commutativity preserving
if and only if $\Phi(A)=ZA$ for all $A\in{\mathcal A}$, where
$Z\in{\mathcal Z}_S({\mathcal A})$ with $Z^2=I$ (Theorem 3.1). It is
clear that Theorem 2.11 above is also an immediate consequence of
this result.

\section{SSCP maps on prime
rings with involution}

In this section, we discuss the  question of characterizing the
strong skew commutativity preserving maps  on prime rings with
involution $*$. The following is our main result in this section.

{\bf Theorem 2.1.} {\it  Let ${\mathcal A}$ be a unital prime *-ring
with the unit $I$. Assume that $\mathcal A$ contains a nontrivial
symmetric idempotent $P$ and $\Phi: {\mathcal A}\rightarrow
{\mathcal A}$ is a surjective map. If $\Phi$ is strong skew
commutativity preserving, that is,
$\Phi(A)\Phi(B)-\Phi(B)\Phi(A)^*=AB-BA^*$ for all $A,B\in{\mathcal
A}$, then there exists a map $f:{\mathcal A}\rightarrow {\mathcal
Z}_S({\mathcal A})$, the set of central symmetric elements,  such
that $\Phi(A)= A+f(A)$ for all $A\in{\mathcal A}$, or $\Phi(A)=
-A+f(A)$ for all $A\in{\mathcal A}$.}

To prove Theorem 2.1, we need a lemma.

Let $\mathcal A$ be a prime ring. Denote by  ${\mathcal Q}={\mathcal
Q}_{ml}({\mathcal A})$ the maximal left ring of quotients of
$\mathcal A$. Note that the center $\mathcal C$ of $\mathcal Q$ is a
field which is called the extended centroid of $\mathcal A$ (see
\cite{BMM, Br} for details). Moreover, ${\mathcal Z}({\mathcal
A})\subseteq{\mathcal C}$. The following result is  well-known.

{\bf Lemma 2.2.} (\cite[Theorem A.7]{Br}) {\it Let $\mathcal A$ be a
prime ring, and let $A_i$, $B_i$, $C_j$, $D_j\in{\mathcal
Q}_{ml}({\mathcal A})$ be such that
$$\sum_{i=1}^n A_iXB_i=\sum_{j=1}^mC_jXD_j\quad{\rm for\ \ all}\ \ X\in{\mathcal A}.$$
If $A_1,\ldots,A_n$ are linearly independent over $\mathcal C$, then
each $B_i$ is a $\mathcal C$-linear combination of $D_1,\ldots,D_m$.
Similarly, if $B_1,\ldots,B_n$ are linearly independent over
$\mathcal C$, then each $A_i$ is a $\mathcal C$-linear combination
of $C_1,\ldots,C_m$. In particular, for $A,B\in{\mathcal
Q}_{ml}({\mathcal A})$, if $AXB=BXA$ for all $X\in{\mathcal A}$,
then $A$ and $B$ are $\mathcal C$-linear dependent. }

\if false The ``if'' part of the theorem is obvious. In the sequel,
we always assume that $\Phi: {\mathcal A}\rightarrow {\mathcal A}$
is a surjective strong skew commutativity preserving map. \fi

 We will prove Theorem 2.1 by a series
of lemmas. Assume in the sequel that $\Phi: {\mathcal
A}\rightarrow{\mathcal A}$ is a SSCP surjective map.

{\bf Lemma 2.3.} {\it $\Phi({\mathcal Z}_S({\mathcal A}))={\mathcal
Z}_S({\mathcal A})$.}

{\bf Proof.} Take any $Z\in {\mathcal Z}_S({\mathcal A})$. Then for
any $T\in{\mathcal A}$, we have $[\Phi(Z),\Phi(T)]_*=[Z,T]_*=0$. So
$\Phi(Z)\Phi(T)=\Phi(T)\Phi(Z)^*.$ By the surjectivity  of $\Phi$,
we get $$\Phi(Z)X=X\Phi(Z)^*\quad{\rm for \ \ all}\quad
X\in{\mathcal A}.\eqno(2.1)$$Take $X=I$ in Eq.(2.1), we get
$\Phi(Z)=\Phi(Z)^*$. This and Eq.(2.1) imply that $\Phi(Z)\in
{\mathcal Z}_S({\mathcal A})$. On the other hand, if there exists
some $A\in{\mathcal A}$ such that  $\Phi(A)=Z=Z^*\in{\mathcal
Z}_S({\mathcal A})$, then we get $[A,S]_*=[\Phi(A),\Phi(S)]_*=0$ for
all $S\in{\mathcal A}$. Hence $AS=SA^*$ holds for all $S$, which
means that $A=A^*\in {\mathcal Z}_S({\mathcal A})$, that is,
${\mathcal Z}_S({\mathcal A})\subseteq\Phi({\mathcal Z}_S({\mathcal
A}))$, completing the proof. \hfill$\Box$

{\bf Lemma 2.4.} {\it For any $T,S\in{\mathcal A}$, there exists an
element $Z_{T,S}\in {\mathcal Z}_S({\mathcal A})$ depending on $T,S$
such that $\Phi(T+S)=\Phi(T)+\Phi(S)+Z_{T,S}$.}

{\bf Proof.} For any $T,S,R\in{\mathcal A}$, we have
$$\begin{array}{rl}
&[\Phi(T+S)-\Phi(T)-\Phi(S),\Phi(R)]_*\\
=&[\Phi(T+S),\Phi(R)]_*-[\Phi(T),\Phi(R)]_*-[\Phi(S),\Phi(R)]_*\\
=&[T+S,R]_*-[T,R]_*-[S,R]_*=0.
\end{array}$$
By the surjectivity of $\Phi$ and the above equation, one sees that
$Z_{T,S}=\Phi(T+S)-\Phi(T)-\Phi(S)\in {\mathcal Z}_S({\mathcal A})$
holds for all $T,S\in{\mathcal A}$.  \hfill$\Box$

In the following, we will use the technique of Peirce decomposition.
By the assumption, we can take   a symmetric nontrivial idempotent
element $P$ in $\mathcal A$. Set ${\mathcal A}_{11}=P{\mathcal A}P$,
${\mathcal A}_{12}=P{\mathcal A}(I-P)$, ${\mathcal
A}_{21}=(I-P){\mathcal A}P$ and ${\mathcal A}_{22}=(I-P){\mathcal
A}(I-P)$. Then ${\mathcal A}={\mathcal A}_{11}+{\mathcal
A}_{12}+{\mathcal A}_{21}+{\mathcal A}_{22}$. It is clear that
${\mathcal A}_{ij}^*={\mathcal A}_{ji}$, $i,j=1,2$. For an element
$S_{ij}\in{\mathcal A}$, we always mean that $S_{ij}\in {\mathcal
A}_{ij}$.

{\bf Lemma 2.5.} {\it $\Phi(P)^*=\Phi(P)$ and
$\Phi(I-P)^*=\Phi(I-P)$. Moreover, there exist elements $\alpha$,
$\beta$, $\mu\in{\mathcal C}$ with $\alpha\not=0$ such that
$\Phi(P)=\alpha P+\mu I$ and $\Phi(I-P)=\alpha(I-P)+\beta I$.}

{\bf Proof.} For any $A\in{\mathcal A}$, it is easy to check that
$[P,[P,[P,A]_*]_*]_*=[P,A]_*$. So we have
$[P,[P,[\Phi(P),\Phi(A)]_*]_*]_*=[\Phi(P),\Phi(A)]_*$. It follows
from the surjectivity of $\Phi$ that
$$[P,[P,[\Phi(P),A]_*]_*]_*=[\Phi(P),A]_*\quad{\rm for\ all}\ A\in{\mathcal
A}.$$ Write $\Phi(P)=S_{11}+S_{12}+S_{21}+S_{22}$. Then the above
equation becomes
$$\begin{array}{rl}&PAS_{11}^*+PAS_{12}^*-PAS_{21}^*-PAS_{22}^*\\
&-S_{11}AP-S_{12}AP+S_{21}AP+S_{22}AP\\
=&S_{21}A+S_{22}A-AS_{21}^*-AS_{22}^*\end{array}\eqno(2.2)$$  for
all $ A\in{\mathcal A}$.

Taking $A=A_{12}$ in Eq.(2.2), we get
$A_{12}S_{12}^*=S_{21}A_{12}=0$, that is,
$$PA(I-P)S_{12}^*=S_{21}PA(I-P)=0\quad{\rm for \ \ all}\ \ A\in{\mathcal A}. $$
It follows from the primeness of $\mathcal A$ that
$S_{12}^*=S_{21}=0$, and so $S_{12}=S_{21}=0$.

Now let $A=A_{11}$ in Eq.(2.2), and we get
$A_{11}S_{11}^*=S_{11}A_{11}$.  This implies that $S_{11}=S_{11}^*$
by taking $A_{11}=P$.  So $PAS_{11}=S_{11}AP$ holds for all
$A\in{\mathcal A}$. It follows from Lemma 2.2 that $S_{11}=\lambda
P$ for some $\lambda\in{\mathcal C}$.

Similarly, taking $A=A_{22}$ in Eq.(2.2), and one can obtain
$S_{22}=S_{22}^*$ and $S_{22}=\mu(I-P)$ for some $\mu\in{\mathcal
C}$. Hence $\Phi(P)=S_{11}+S_{22}=S_{11}^*+S_{22}^*=\Phi(P)^*$ and
$$\Phi(P)=S_{11}+S_{22}=\lambda P+\mu(I-P)=\alpha P+\mu I,$$
where $\alpha=\lambda-\mu\in{\mathcal C}$. It is obvious that $\mu
I\in{\mathcal C}$. Note that $\Phi(I)\in{\mathcal Z}_S({\mathcal
A})$, $\Phi(I)-\Phi(P)-\Phi(I-P)\in{\mathcal Z}_S({\mathcal A})$ and
${\mathcal Z}_S({\mathcal A})\subseteq{\mathcal C}$. So
$\Phi(I-P)^*=\Phi(I-P)$ and there exists an element
$\beta\in{\mathcal C}$ such that
$$\Phi(I-P)=\alpha(I-P)+\beta I.$$

Finally, we still need to prove that $\alpha\not=0$. On the
contrary, if $\alpha=0$, then $\Phi(P)=\mu I\in{\mathcal C}$. Since
$\Phi(P)\in{\mathcal A}$, it follows that $\Phi(P)\in{\mathcal
Z}_S({\mathcal A})$. By Lemma 2.2, we get $P\in{\mathcal
Z}_S{(\mathcal A)}$, which is impossible as $\mathcal A$ is prime.
The proof is finished. \hfill$\Box$

Note that $\mathcal C$ is a field as $\mathcal A$ is prime
(\cite[Theorem A.6]{Br}). So $\alpha\in{\mathcal C}$ is invertible.
In the following, let $\lambda=\alpha^{-1}\in{\mathcal C}$. Also
note that the unit $1$ of ${\mathcal C}$ is the same to the unit $I$
of ${\mathcal A}$

{\bf Lemma 2.6.} {\it For any $A_{ij}\in{\mathcal A}_{ij}$, we have
$\Phi(A_{ij})=\lambda A_{ij}$, $1\leq i\not=j\leq 2$.}

{\bf Proof.} Take any $A_{12}\in{\mathcal A}_{12}$ and let
$\Phi(A_{12})=S_{11}+S_{12}+S_{21}+S_{22}$. Since
$[\Phi(P),\Phi(A_{12})]_*=[P,A_{12}]_*=A_{12}$, by Lemma 2.5, we get
$\alpha S_{12}-\alpha S_{21}=A_{12}$, which implies that $S_{21}=0$
and $S_{12}=\alpha^{-1}A_{12}=\lambda A_{12}$.

For any $B\in{\mathcal A}$, by the surjectivity of $\Phi$, there
exists an element $T=T_{11}+T_{12}+T_{21}+T_{22}\in{\mathcal A}$
such that $\Phi(T)=B$. Since
$[B,\Phi(A_{12})]_*=[\Phi(T),\Phi(A_{12})]_*=[T,A_{12}]_*$, we have
$$\begin{array}{rl}&BS_{11}+ B(\lambda A_{12})+BS_{22}-S_{11}B^*-\lambda A_{12}B^*-S_{22}B^*\\
=&T_{11}A_{12}+T_{21}A_{12}-A_{12}T_{12}^*-A_{12}T_{22}^*.\end{array}\eqno(2.3)$$
Multiplying by $P$ from the right in Eq.(2.3), one gets
$$BS_{11}-S_{11}B^*P-\lambda A_{12}B^*P-S_{22}B^*P
=-A_{12}T_{12}^*.$$ Replacing $B$ by $(I-P)BP$ in the above
equation, we have $(I-P)BPS_{11} =-A_{12}T_{12}^*$, which implies
that $(I-P)BPS_{11}=0$ for all $B\in{\mathcal A}$. It follows from
the primeness of ${\mathcal A}$ that $S_{11}=0$.

Similarly, replacing $B$ by $PB(I-P)$  and multiplying by $I-P$ from
the left in Eq.(2.3), one can show that $S_{22}=0$. Hence we obtain
$\Phi(A_{12})=S_{12}=\lambda A_{12}$.

The proof of  $\Phi(A_{21})=\lambda A_{21}$ is similar, and we omit
it here.  \hfill$\Box$

{\bf Lemma 2.7.} {\it For any $A_{ii}\in{\mathcal A}_{ii}$, we have
$\Phi(A_{ii})=\alpha A_{ii}$, $i=1,2$.}

{\bf Proof.} Still, we only need to prove that the lemma is true for
$A_{11}$.

Take any $A_{11}\in{\mathcal A}_{11}$ and let
$\Phi(A_{11})=S_{11}+S_{12}+S_{21}+S_{22}$. Since
$[\Phi(P),\Phi(A_{11})]_*=[P,A_{11}]_*=0$, by Lemma 2.5, we have
$$\alpha S_{11}+\alpha S_{12}-S_{11}(\alpha P)-S_{21}(\alpha P)=0.$$
This implies $\alpha S_{12}=0$, an so $S_{12}=0$. On the other hand,
since $[\Phi(I-P),\Phi(A_{11})]_*=[I-P,A_{11}]_*=0$, by Lemma 2.5
again, one gets $\alpha S_{21}=0$, which implies that $S_{21}=0$.

For  any $B_{12}\in{\mathcal A}_{12}$, by Lemma 2.6, we have
$[\lambda
B_{12},\Phi(A_{11})]_*=[\Phi(B_{12}),\Phi(A_{11})]_*=[B_{12},A_{11}]_*=0,$
that is, $$\lambda B_{12}S_{22}-S_{11}(\lambda
B_{12})^*-S_{22}(\lambda B_{12})^*=0.$$ Note that $\lambda
B_{12}\in{\mathcal A}_{12}$. It follows that $\lambda
B_{12}S_{22}=0$, which implies  $B_{12}S_{22}=0$. That is,
$PB(I-P)S_{22}=0$ for all $B\in{\mathcal A}$. As ${\mathcal A}$ is
prime, we get $S_{22}=0$.

For  any $B_{21}\in{\mathcal A}_{21}$, by Lemma 2.6, we have
$[\lambda
B_{21},\Phi(A_{11})]_*=[\Phi(B_{21}),\Phi(A_{11})]_*=[B_{21},A_{11}]_*$,
that is,
$$\lambda B_{21}S_{11}-S_{11}(\lambda B_{21})_*=B_{21}A_{11}-A_{11}B_{21}^*.$$
This forces $B_{21}(\lambda S_{11}-A_{11})=(\lambda
S_{11}-A_{11})B_{21}^*=0$. So we get $(I-P)BP(\lambda
S_{11}-A_{11})=0$ for each $B\in{\mathcal A}$. It follows from the
primeness of ${\mathcal A}$ that $\lambda S_{11}=A_{11}$. Hence
$\Phi(A_{11})=S_{11}=\lambda^{-1}A_{11}=\alpha A_{11}$, completing
the proof. \hfill$\Box$

{\bf Lemma 2.8.} {\it $\alpha=\lambda$, and hence $\alpha^2=1$ and
$\alpha =\pm 1$.}

{\bf Proof.} For any $A_{12}\in {\mathcal A}_{12}$ and $A_{21}\in
{\mathcal A}_{21}$, by the definition of $\Phi$ and Lemma 2.6, we
have
$$A_{12}A_{21}=[A_{12}, A_{21}]_*=[\Phi(A_{12}),\Phi(A_{21})]_*=[\lambda
A_{12},\lambda A_{21}]_*=(\lambda A_{12})(\lambda
A_{21})=\lambda^2A_{12} A_{21}.$$It follows that
$(\lambda^2-1)A_{12}A_{21}=0$. First fix $A_{12}$. Then the
 equation becomes $(\lambda^2-1)A_{12}AP=0$ for all $A\in{\mathcal
A}$. Assume that $\lambda^2-1\not=0$. Since $\mathcal C$ is a field,
we get that $\lambda^2-1$ is invertible. So we have $A_{12}AP=0$ for
all $A\in{\mathcal A}$. Since $\mathcal A$ is prime, it follows that
$A_{12}=0$, that is, $PA(I-P)=0$ for all $A\in{\mathcal A}$. This
implies that either $P=0$ or $P=I$, a contradiction. So
$\lambda^2=1$ and $\lambda=\alpha$. Since $\mathcal C$ is a field,
we see that $\alpha =1$ or $-1$, completing the proof. \hfill$\Box$

{\bf Lemma 2.9.} {\it There exists a map $f: {\mathcal A}\rightarrow
{\mathcal Z}({\mathcal A})$   such that $\Phi(A)=A+f(A)$ for all
$A\in{\mathcal A}$ or $\Phi(A)=-A+f(A)$ for all $A\in{\mathcal A}$.}

{\bf Proof.} By Lemmas 2.6-2.8, we have proved that
$\Phi(A_{ij})=A_{ij}$  for all $A_{ij}\in{\mathcal A}_{ij}$ or
$\Phi(A_{ij})=-A_{ij}$  for all $A_{ij}\in{\mathcal A}_{ij}$
($i,j=1,2$). Now, for any $A=A_{11}+A_{12}+A_{21}+A_{22}\in
{\mathcal A}$, by Lemma 2.4, there exists some element
$Z_A\in{\mathcal Z}_S({\mathcal A})$ such that
$$\begin{array}{rl}
&\Phi(A)-(\Phi(A_{11})+\Phi(A_{12})+\Phi(A_{21})+\Phi(A_{22}))\\
=&\Phi(A_{11}+A_{12}+A_{21}+A_{22})-\Phi(A_{11})-\Phi(A_{12})-\Phi(A_{21})-\Phi(A_{22})=Z_A.
\end{array}$$
Define a map $f: {\mathcal A}\rightarrow {\mathcal Z}_S({\mathcal
A})$ by $f(A)=Z_A$ for each $A\in{\mathcal A}$. Then
$\Phi(A)=A+f(A)$ for each $A$, or $\Phi(A)=-A+f(A)$ for each $A$.
\hfill$\Box$

 By applying Theorem 2.1, we give  a characterization of SSCP
 maps
 on prime *-algebras with the second kind involution.

{\bf Theorem 2.10.} {\it  Let ${\mathcal A}$ be a unital prime
*-algebra over a field $\mathbb F$ with a nontrivial symmetric idempotent $P$. Assume that the involution $*$ is of the second kind and $\Phi:
{\mathcal A}\rightarrow {\mathcal A}$ is a surjective map. Then
$\Phi$ is strong skew  commutativity  preserving (that is,
$\Phi(A)\Phi(B)-\Phi(B)\Phi(A)^*=AB-BA^*$ for all $A,B\in{\mathcal
A}$) if and only if  $\Phi(A)=A$ for all $A\in{\mathcal A}$ or
$\Phi(A)=-A$ for all $A\in{\mathcal A}$.}

{\bf Proof.} Obviously, the ``if'' part is true. For the ``only if''
part, by Theorem 2.1, there exists a map $f: {\mathcal A}\rightarrow
{\mathcal Z}_S({\mathcal A})$ such that $\Phi(A)=A+f(A)$ for all
$A\in{\mathcal A}$, or  $\Phi(A)=-A+f(A)$ for all $A\in{\mathcal
A}$. So, to complete the proof of the theorem, we only need to prove
$f\equiv0$.

In fact, $\Phi$ is SSCP implies that
$$f(B)(A-A^*)=0\quad{\rm for\ \ all}\quad A,B\in{\mathcal A}.\eqno(2.4)$$
Since $*$ is of the second kind, there exists a nonzero
$\epsilon\in{\mathbb F}$ such that $(\epsilon I)^*=-\epsilon I$.
Thus, let $A=\epsilon I$ in Eq.(2.4), we get $2\epsilon f(B)=0$,
which implies $f(B)=0$ for each $B\in{\mathcal A}$ as $\mathbb F$ is
a field. The proof is finished. \hfill$\Box$

Recall that a von Neumann algebra ${\mathcal M}$  is called a factor
if its center is trivial (i.e., ${\mathcal Z}({\mathcal M})={\mathbb
C}I$). Note that factor von Neumann algebras are prime and $*$ is of
the second kind. So, as an application of Theorem 2.10 to the factor
von Neumann algebras case, we improve the main results of \cite{CP}
immediately.

\textbf{Theorem 2.11.} {\it Let ${\mathcal A}$ be a factor von
Neumann algebra. Assume that $\Phi: {\mathcal A}\rightarrow{\mathcal
A}$ is a surjective map. Then $\Phi$ is strong skew commutativity
preserving if and only if $\Phi(A)=\alpha A$ for all $A\in{\mathcal
A}$, where $\alpha\in \{1, -1\}$.}

\section{SSCP maps on von Neumann
algebras}

In this section, we will discuss the  strong skew commutativity
preserving maps on von Neumann algebras. The following is the main
result of this section.

{\bf Theorem 3.1.} {\it  Let ${\mathcal M}$ be a  von Neumann
algebra without central summands of type $I_1$. Assume that $\Phi:
{\mathcal M}\rightarrow {\mathcal M}$ is a surjective map. Then
$\Phi$ is strong skew commutativity preserving if and only if there
exists an element $Z\in{\mathcal Z}_S({\mathcal M})$ with $Z^2=I$
such that $\Phi(A)=ZA$ for all $A\in{\mathcal A}$.}

Obviously, Theorem 2.11 is also an immediate consequence of the
above result.

We remark that the methods used in the proof of Theorem 2.1 are not
valid here since the von Neumann algebras in Theorem 3.1 may not be
prime. In order to overcome the difficulties caused by the absence
of primeness, we need some deep results coming from the theory of
von Neumann algebras.

Let $\mathcal M$ be a von Neumann algebra and $A\in{\mathcal M}$.
Recall that the central carrier of $A$, denoted by $\overline{A}$,
is the intersection of all central projections $P$ such that $PA=A$.
If $A$ is self-adjoint, then the core of $A$, denoted by
$\underline{A}$, is sup$\{S\in{\mathcal Z}({\mathcal M}):S=S^*,
S\leq A\}$. Particularly, if $A=P$ is a projection, it is clear that
$\underline{P}$ is the largest central projection $\leq P$. A
projection $P$ is said to be core-free if $\underline{P}=0$
\cite{Mi}. It is easy to see that $\underline{P}=0$ if and only if
$\overline{I-P}=I$.

The following lemmas are useful for our purpose, where Lemma 3.2 and
Lemma 3.3 are proved in \cite{Mi}.

{\bf Lemma 3.2.} (\cite{Mi}) {\it Let $\mathcal M$ be a von Neumann
algebra without central summands of type $I_1$. Then each nonzero
central projection $C\in{\mathcal M}$ is the carrier of a core-free
projection in $\mathcal M$. Particularly, there exists a nonzero
core-free projection $P\in{\mathcal M}$ with $\overline{P}=I$.}

{\bf Lemma 3.3.} (\cite{Mi}) {\it Let $\mathcal M$ be a von Neumann
algebra. For projections $P,Q\in{\mathcal M}$, if
$\overline{P}=\overline{Q}\not=0$ and $P+Q=I$, then $T\in{\mathcal
M}$ commutes with $PXQ$ and $QXP$ for all $X\in{\mathcal M}$ implies
that $T\in{\mathcal Z}({\mathcal M})$.}

{\bf Lemma 3.4.} {\it Let $\mathcal M$ be a  von Neumann algebra.
Assume that $P\in\mathcal M$ is a projection satisfying
$\underline{P}=0$ and $\overline{P}=I$.  Then, for any
$Z\in{\mathcal Z}({\mathcal M})$,  $ZP{\mathcal M}(I-P)=\{0\}$
implies $Z=0$.}

{\bf Proof.} Assume that ${\mathcal M}\subseteq{\mathcal B}(H)$,
where $H$ is a Hilbert space, and assume that $Z\in{\mathcal
Z}({\mathcal M})$ with $Z\not=0$ such that $ZP{\mathcal
M}(I-P)=\{0\}$. Let $Q$ be the projection onto the closure of the
range of $Z$. It is clear that $Q\in{\mathcal Z}({\mathcal M})$. So
we may write ${\mathcal M}=Q{\mathcal M}Q\oplus (I-Q){\mathcal
M}(I-Q)={\mathcal M}_1\oplus {\mathcal M}_2$. Thus, according to the
space decomposition $H=QH\oplus (I-Q)H$, we have
$$A=\left(\begin{array}{rl}A_1&0\\0&A_2\end{array}\right),\quad
Z=\left(\begin{array}{rl}Z_1&0\\0&0\end{array}\right),\quad
P=\left(\begin{array}{rl}P_1&0\\0&P_2\end{array}\right),$$ where
 $A\in{\mathcal M}$ is arbitrary, $A_i,P_i\in{\mathcal M}_i$ with $P_i=P_i^*=P_i^2$ ($i\in\{1,2\}$)
 and $Z_1\in{\mathcal Z}({\mathcal M}_1)$ is   injective with dense range.
It follows that
$$ZPA(I-P)=\left(\begin{array}{cc}Z_1P_1A_1(I_1-P_1)&0\\0&0\end{array}\right).$$
By the assumption, we get $Z_1P_1A_1(I_1-P_1)=0$, which implies that
$P_1A_1(I_1-P_1)=0$ for all $A_1\in{\mathcal M}_1$ as $Z_1$ is
injective. So $P_1A_1=P_1A_1P_1$ for each $A_1\in{\mathcal M}_1$.
Thus, by the following claim, we have $P_1\in{\mathcal Z}({\mathcal
M}_1)$.

{\bf Claim.} Let $\mathcal N$ be a  von Neumann algebra.  Assume
that $P\in{\mathcal N}$ is a projection satisfying $P{\mathcal
N}(I-P)=\{0\}$. Then $P\in{\mathcal Z}({\mathcal N})$.

In fact, for such $P$, write $P_1=P$ and $P_2=I-P$. Denote
${\mathcal N}_{ij}=P_i{\mathcal N}P_j$, $i,j\in\{1,2\}$. Then
${\mathcal N}={\mathcal N}_{11}+{\mathcal N}_{12}+{\mathcal
N}_{21}+{\mathcal N}_{22}$. For any $A\in{\mathcal N}$, we have
$PA=PAP+PA(I-P)$. It follows from the assumption $P{\mathcal
N}(I-P)=\{0\}$  that $PA(I-P)=0$ for any $A\in{\mathcal N}$. Since
$A$ is arbitrary, it is true that $PA^*(I-P)=0$ holds for any
$A\in{\mathcal N}$, which implies that $(I-P)AP=0$ for any
$A\in{\mathcal N}$. Thus, for any $A\in{\mathcal N}$, we must have
$A=PAP+(I-P)A(I-P)$. Now it is clear that $PA=AP$ for each $A$, that
is, the claim is true.

Let us go back to the proof of the lemma and let
$Q_0=\left(\begin{array}{rl}P_1&0\\0&0\end{array}\right)$.
Obviously, $Q_0$ is a central projection  with $Q_0\leq P$. Since
$\underline{P}=0$, we have $Q_0=0$, and so $P_1=0$. This yields
$QP=0$, which implies $I-Q\geq \overline{P}=I$, a contradiction.
Hence $Z=0$ and the proof is completed. \hfill$\Box$

Now we are in a position to give our proof of Theorem 3.1.

{\bf Proof of Theorem 3.1.} Still, we only need to prove the ``only
if'' part.

By the same argument as that of Lemmas 2.3-2.4, one can obtain that
$$\Phi({\mathcal Z}_S({\mathcal M}))={\mathcal Z}_S({\mathcal M})\eqno(3.1)$$
and$$\Phi(T+S)=\Phi(T)+\Phi(S)+Z_{T,S}\quad{\rm for \ \ all}\quad
T,S\in{\mathcal M},\eqno(3.2)$$ where  $Z_{T,S}\in {\mathcal
Z}_S({\mathcal M})$ depending on $T,S$.

By Lemma 3.2, there is a non-central core-free projection $P$ with
central carrier $I$. For such a $P$,  by the definitions of core and
central carrier, $I-P$ is also core-free with $\overline{I-P}=I$.
Denote ${\mathcal M}_{ij}=P_i{\mathcal M}P_j$ ($i,j\in\{1,2\}$),
where $P_1=P$ and $P_2=I-P$. Then ${\mathcal M}={\mathcal
M}_{11}+{\mathcal M}_{12}+{\mathcal M}_{21}+{\mathcal M}_{22}$. In
all that follows, when  writing $S_{ij}$,  it always indicates
$S_{ij}\in {\mathcal M}_{ij}$.

We will complete the proof by   several steps.

{\bf Step 1.} $\Phi(P)=\Phi(P)^*$.

For the identity operator $I$, by Eq.(3.1), there exists some
$Z\in{\mathcal Z}_S({\mathcal M})$ such that $\Phi(Z)=I$. Thus
$$0=[P,Z]_*=[\Phi(P),\Phi(Z)]_*=\Phi(P)-\Phi(P)^*,$$
which implies $\Phi(P)=\Phi(P)^*$.

{\bf Step 2.} There exist elements $Z_1,Z_2,Z_3\in{\mathcal
Z}_S({\mathcal M})$ with $Z_1\not=0$ such that $\Phi(P)=Z_1P+Z_2$
and $\Phi(I-P)=Z_1(I-P)+Z_3$.

Note that $[P,[P,[P,A]_*]_*]_*=[P,A]_*$ holds for any $A\in{\mathcal
M}$. Thus for every $A$, we have
$$[P,[P,[\Phi(P),\Phi(A)]_*]_*]_*=[\Phi(P),\Phi(A)]_*.$$
It follows from the surjectivity of $\Phi$ that
$$[P,[P,[\Phi(P),A]_*]_*]_*=[\Phi(P),A]_*\quad{\rm for\ all}\ A\in{\mathcal
M}.\eqno(3.3)$$ Write $[\Phi(P),A]_*=B$. Then Eq.(3.3) becomes
$$PB-2PBP+BP=B.$$
Multiplying by $P$ and $I-P$ from both sides in the above equation,
respectively,  one gets $PBP=0$ and $(I-P)B(I-P)=0$. Therefore
$$P(\Phi(P)A-A\Phi(P)^*)P=0\quad{\rm and}\quad (I-P)(\Phi(P)A-A\Phi(P)^*)(I-P)=0\eqno(3.4)$$
hold for all $A\in{\mathcal M}$.

Write $\Phi(P)=S_{11}+S_{12}+S_{21}+S_{22}$. Replacing $A$ by
$PT(I-P)$ for any $T\in{\mathcal M}$ in Eq.(3.4), we get
$$PT(I-P)S_{12}^*=0\quad{\rm and}\quad S_{21}PT(I-P)=0,$$
which, together with $\Phi(P)=\Phi(P)^*$, implies  that
$$PT(I-P)S_{21}=S_{21}PT(I-P)=0\eqno(3.5)$$
holds for all $T\in{\mathcal M}$.  It is obvious that
$S_{21}(I-P)TP=(I-P)TPS_{21}=0$. Hence, by Lemma 3.3, we see that
$S_{21}\in{\mathcal Z}({\mathcal M})$, which forces $S_{21}=0$.

Similarly, replacing $A$ by $(I-P)TP$ for any $T\in{\mathcal M}$ in
Eq.(3.4) and applying Step 1,  it is easily checked $S_{12}=0$.

Now taking $A=P$ in Eq.(3.3), and by Step 1, one gets
$$S_{11}PTP=PTPS_{11}^*=PTPS_{11},$$ which means
$S_{11}\in P{\mathcal Z}_S({\mathcal M})$. Symmetrically, by taking
$A=I-P$ in Eq.(3.3), we obtain $S_{22}\in (I-P){\mathcal
Z}_S({\mathcal M})$. Write $S_{11}=Z_{11}P$ and
$S_{22}=Z_{22}(I-P)$, where $Z_{11},Z_{22}\in{\mathcal
Z}_S({\mathcal M})$. Then
$$\Phi(P)=S_{11}+S_{22}
=Z_{11}P+Z_{22}(I-P)=(Z_{11}-Z_{22})P+Z_{22}.$$ Note that
$\Phi(I)-\Phi(I-P)-\Phi(P)\in{\mathcal Z}_S({\mathcal M})$ by
Eqs.(3.1)-(3.2). So there exists some $Z_{33}\in{\mathcal
Z}_S({\mathcal M})$ such that
$$\begin{array}{rl}\Phi(I-P)=&Z_{33}-\Phi(P)=Z_{33}-Z_{22}-(Z_{11}-Z_{22})P\\
=&(Z_{11}-Z_{22})(I-P)+(Z_{33}-Z_{11}).\end{array}$$ Let
$Z_1=Z_{11}-Z_{22}$, $Z_2=Z_{22}$ and $Z_3=Z_{33}-Z_{11}$. Thus
$\Phi(P)=Z_1P+Z_2$ and $\Phi(I-P)=Z_1(I-P)+Z_3$.

Finally, we still need to prove that $Z_1\not=0$. On the contrary,
if $Z_1=0$, then $\Phi(P)=Z_2\in{\mathcal Z}({\mathcal M})$.  By
Eq.(3.1), we get $P\in{\mathcal Z}_S{(\mathcal M)}$, which is
impossible.

{\bf Step 3.} If $\Phi(T)=P$ and $\Phi(S)=I-P$, where
$T=T_{11}+T_{12}+T_{21}+T_{22}\in{\mathcal M}$ and
$S=S_{11}+S_{12}+S_{21}+S_{22}\in{\mathcal M}$, then
$T_{12}=T_{21}=0$ and $S_{12}=S_{21}=0$.

\if false Since $[T,P]_*=[\Phi(T),\Phi(P)]_*=[P,\Phi(P)]_*$, by Step
2, we have $T_{11}+T_{21}-T_{11}^*-T_{21}^*=0$, which implies
$T_{21}=0$ and $T_{11}=T_{11}^*$.\fi

In fact, by the equation $[P,T]_*=[\Phi(P),\Phi(T)]_*=[\Phi(P),P]_*$
and Step 2, one can get $T_{12}=T_{21}=0$;
 by the equation $[I-P,S]_*=[\Phi(I-P),\Phi(S)]_*=[\Phi(I-P),I-P]_*$
and Step 2, one can get $S_{12}=S_{21}=0$.

{\bf Step 4.}  For any $A_{ij}\in{\mathcal M}_{ij}$,  we have
$\Phi(A_{ij})\in{\mathcal M}_{ij}$, $1\leq i\not=j\leq 2$. Moreover,
$Z_1P\Phi(A_{12})(I-P)=A_{12}$ holds for all $A_{12}\in{\mathcal
M}_{12} $ and $Z_1(I-P)\Phi(A_{21})P=A_{21}$ holds for all
$A_{21}\in{\mathcal M}_{21} $.

We only check the assertion for $A_{12}$, and the case of $A_{21}$
is similarly dealt with.

For any $A_{12}$, write $\Phi(A_{12})=S_{11}+S_{12}+S_{21}+S_{22}$.
By Step 3, there exists some $T=T_{11}+T_{22}\in{\mathcal M}$ such
that $\Phi(T)=P$. Then $[P,\Phi(A_{12})]_*=[T,A_{12}]_*$, that is,
$$S_{12}-S_{21}=T_{11}A_{12}-A_{12}T_{22}^*.$$
Multiplying by $I-P$ and $P$ from the left side and the right side
respectively  in the above equation, one gets
$$S_{21}=(I-P)\Phi(A_{12})P=0.\eqno(3.6)$$

On the other hand, since $[\Phi(A_{12}),P]_*=[A_{12},T]_*$, by
Eq.(3.6), we have $S_{11}-S_{11}^*=A_{12}T_{22}-T_{22}A_{12}^*.$
Multiplying by $P$ from both sides of the equation, it follows that
$$S_{11}=S_{11}^*.\eqno(3.7)$$

Similarly,  by using the equation
$[\Phi(A_{12}),I-P]_*=[\Phi(A_{12}),\Phi(S)]_*=[A_{12},S]_*$ and
Step 3, it is easily checked that
$$S_{22}=S_{22}^*.\eqno(3.8)$$

Now, for  any $X\in{\mathcal M}$, by the surjectivity of $\Phi$,
there exists an element $R=R_{11}+R_{12}+R_{21}+R_{22}\in{\mathcal
M}$ such that $\Phi(R)=X$. Since
$[X,\Phi(A_{12})]_*=[\Phi(R),\Phi(A_{12})]_*=[R,A_{12}]_*$, applying
Eq.(3.6), we get
$$\begin{array}{rl}&XS_{11}+XS_{12}+XS_{22}-S_{11}X^*-S_{12}X^*-S_{22}X^*\\
=&R_{11}A_{12}+R_{21}A_{12}-A_{12}R_{12}^*-A_{12}R_{22}^*.\end{array}\eqno(3.9)$$
Replacing $X$ by $PB(I-P)+(I-P)BP$  for any $B\in{\mathcal M}$, and
multiplying by $I-P$ and $P$ from the left and the right
respectively in Eq.(3.9), one obtains
$$(I-P)BPS_{11}=S_{22}(I-P)B^*P\quad\mbox{\rm holds for all }\quad B\in{\mathcal M};\eqno(3.10)$$
Replacing $X$ by $(I-P)BP$  for any $B\in{\mathcal M}$, and
multiplying by $I-P$ and $P$ from the left and the right
respectively in Eq.(3.9),  one gets
$$(I-P)BPS_{11}=0\quad{\rm holds\ for  \ all}\quad B\in{\mathcal M}.\eqno(3.11)$$
Equivalently,
$$(I-P)B^*PS_{11}=0\quad{\rm holds\ for \  all}\quad B\in{\mathcal M}.\eqno(3.12)$$
By Eqs.(3.7) and (3.12), we have
$$S_{11}PB(I-P)=S_{11}^*PB(I-P)=((I-P)B^*PS_{11})^*=0.\eqno(3.13)$$
Moreover, it is obvious that $PB(I-P)S_{11}=S_{11}(I-P)BP=0$. This,
together with Eqs.(3.11), (3.13) and Lemma 3.3, yields
$S_{11}\in{\mathcal Z}_S({\mathcal M})$. Hence
$$S_{11}=0.\eqno(3.14)$$

Combining Eq.(3.10) with Eq.(3.11), we get
$$S_{22}(I-P)B^*P=0\quad{\rm for \ \ all}\quad B\in{\mathcal M}.\eqno(3.15)$$
Equivalently,
$$S_{22}(I-P)BP=0\quad{\rm for \ \ all}\quad B\in{\mathcal M}.$$
It follows from Eqs.(3.8) and (3.15) that
$$PB(I-P)S_{22}=PB(I-P)S_{22}^*
=(S_{22}(I-P)B^*P)^*=0.$$ Thus we have
$$\begin{array}{rl}&PB(I-P)S_{22}=S_{22}PB(I-P)\\
=&S_{22}(I-P)BP=(I-P)BP(I-P)S_{22}=0,\end{array}$$ which implies
that $S_{22}\in{\mathcal Z}_S({\mathcal M})$ by Lemma 3.3, and so
$$S_{22}=0.\eqno(3.16)$$
Combining Eqs.(3.6), (3.14) and (3.16), we see that
$\Phi(A_{12})\in{\mathcal M}_{12}$, as desired.

Finally, since $[\Phi(P),\Phi(A_{12})]_*=[P,A_{12}]_*$, by Step 2,
it is easy to check that  $A_{12}=Z_1S_{12}=Z_1P\Phi(A_{12})(I-P)$.

{\bf Step 5.}  For any $A_{ii}\in{\mathcal M}_{ii}$, we have
$\Phi(A_{ii})=Z_1A_{ii}$, $i=1,2$.

Still, we only  prove  that $\Phi(A_{11})=Z_1A_{11}$ holds for all
$A_{11}\in{\mathcal A}_{11}$. The case of $i=2$ is checked
similarly.

Take any $A_{11}\in{\mathcal M}_{11}$ and any $B\in{\mathcal M}$.
Write $\Phi(A_{11})=S_{11}+S_{12}+S_{21}+S_{22}$. On the one hand,
since $0=[PB(I-P),A_{11}]_*=[\Phi(PB(I-P)),\Phi(A_{11})]_*$, by Step
4, we have
$$\begin{array}{rl}&P\Phi(PB(I-P))(I-P)S_{21}+P\Phi(PB(I-P))(I-P)S_{22}\\
&-S_{12}(I-P)\Phi(PB(I-P))^*P-S_{22}(I-P)\Phi(PB(I-P))^*P=0,\end{array}$$
which implies that
$$P\Phi(PB(I-P))(I-P)S_{21}=S_{12}(I-P)\Phi(PB(I-P))^*P,$$
$$P\Phi(PB(I-P))(I-P)S_{22}=0$$and $$S_{22}(I-P)\Phi(PB(I-P))^*P=0.$$
Multiplying by $Z_1\in{\mathcal Z}_S({\mathcal M})$  in the above
three equations, and applying Step 4, one gets
$$PB(I-P)(I-P)S_{21}=S_{12}(I-P)B^*P,\eqno(3.17)$$
$$PB(I-P)S_{22}=0\eqno(3.18)$$
and
$$S_{22}(I-P)B^*P=0\eqno(3.19)$$
hold for all $B\in{\mathcal M}$. Furthermore, Eq.(3.19) yields
$$S_{22}(I-P)BP=0\quad{\rm for \ \ all} \quad B\in{\mathcal M}.\eqno(3.20)$$
Note that $(I-P)BP(I-P)S_{22}=0$ and $S_{22}PB(I-P)=0$. These,
together with Lemma 3.3, Eqs.(3.18) and (3.20), imply
$S_{22}\in{\mathcal Z}({\mathcal M})$. So we must have
$$S_{22}=(I-P)\Phi(A_{11})(I-P)=0.\eqno(3.21)$$

On the other hand, by using the equation
$[A_{11},PB(I-P)]_*=[\Phi(A_{11}),\Phi(PB(I-P))]_*,$ we get
$$\begin{array}{rl}A_{11}PB(I-P)=&S_{11}P\Phi(PB(I-P))(I-P)\\&
+S_{21}P\Phi(PB(I-P))(I-P)\\&-P\Phi(PB(I-P))(I-P)S_{12}^*.\end{array}$$
It follows that
$$A_{11}PB(I-P)=S_{11}P\Phi(PB(I-P))(I-P),$$
$$S_{21}P\Phi(PB(I-P))(I-P)=0$$
and
$$P\Phi(PB(I-P))(I-P)S_{12}^*=0.$$
Multiplying by $Z_1$ in the above three equations leads to
$$(Z_1A_{11}-S_{11})PB(I-P)=0,\eqno(3.22)$$
$$S_{21}PB(I-P)=0\eqno(3.23)$$
and
$$PB(I-P)S_{12}^*=0\eqno(3.24)$$
for all $B\in{\mathcal M}$. Eq.(3.24) implies that
$$S_{12}(I-P)B^*P=0,\eqno(3.25)$$ and so
$$S_{12}(I-P)BP=0\quad{\rm for \ \ all}\quad B\in{\mathcal M}.\eqno(3.26)$$
  Eq.(3.17) and Eq.(3.25) together yield
$$PB(I-P)S_{21}=0\quad{\rm for \ \ all}\quad B\in{\mathcal M}.\eqno(3.27)$$
It is obvious that $(I-P)BP(I-P)S_{21}=S_{21}(I-P)BP=0$. By
Eqs.(3.23), (3.27) we get $S_{21}\in{\mathcal Z}({\mathcal M})$ and
hence
$$S_{21}=(I-P)\Phi(A_{11})P=0.\eqno(3.28)$$
Now consider the equation
$[(I-P)BP,A_{11}]_*=[\Phi((I-P)BP),\Phi(A_{11})]_*.$ It follows from
Step 4, Eqs.(3.21) and (3.28) that
$$\begin{array}{rl}&(I-P)BPA_{11}-A_{11}PB^*(I-P)=(I-P)\Phi((I-P)BP)PS_{11}\\&
+(I-P)\Phi((I-P)BP)PS_{12}
-S_{11}P\Phi(PB(I-P))^*(I-P).\end{array}$$ This implies that
$$(I-P)BPA_{11}=(I-P)\Phi((I-P)BP)PS_{11}$$
and
$$(I-P)\Phi((I-P)BP)PS_{12}=0.$$
Multiplying by $Z_1$ in the above two equations, we obtain
$$(I-P)BP(Z_1A_{11}-S_{11})=0\eqno(3.29)$$
and
$$(I-P)BPS_{12}=0.\eqno(3.30)$$
Since $S_{12}PB(I-P)=PB(I-P)PS_{12}=0$, by Eqs.(3.26) and (3.30), we
get $S_{12}\in{\mathcal Z}({\mathcal M})$, and so
$$S_{12}=P\Phi(A_{11})(I-P)=0.\eqno(3.31)$$
By Eqs.(3.22) and (3.29), also noting that
$(Z_1A_{11}-S_{11})(I-P)BP=PB(I-P)(Z_1A_{11}-S_{11})=0$, we see that
$Z_1A_{11}-S_{11}\in{\mathcal Z}({\mathcal M})$, and hence
$$Z_1A_{11}=S_{11}=P\Phi(A_{11})P.\eqno(3.32)$$

Now it is clear by Eqs.(3.21), (3.28), (3.31) and (3.32)   that
$\Phi(A_{11})=Z_1A_{11}$.

{\bf Step 6.}  $Z_1^2=I$.

For any $A,B\in {\mathcal M}$, by the assumptions on $\Phi$ and Step
4, we have
$$\begin{array}{rl}PA(I-P)BP=&[PA(I-P),(I-P)BP]_*\\=&[\Phi(PA(I-P)),\Phi((I-P)BP)]_*\\
=&P\Phi(PA(I-P))(I-P)\Phi((I-P)BP)P.\end{array}$$ Multiplying by
$Z_1^2$ in the above equation and applying  Step 4 again, one gets
$Z^2_1PA(I-P)(I-P)BP=PA(I-P)(I-P)BP$. Let us fix $A$. Then the
equation becomes
$$(Z^2_1-I)PA(I-P)(I-P)BP=0\quad{\rm for\ \ all}\quad B\in{\mathcal M}.\eqno(3.33)$$
Similarly, by using the equation
$[(I-P)BP,PA(I-P)]_*=[\Phi((I-P)BP),\Phi(PA(I-P))]_*$, one   obtains
$$(I-P)BP(Z^2_1-I)PA(I-P)=(Z^2_1-I)(I-P)BPPA(I-P)=0\quad{\rm for\ \ all}\quad B\in{\mathcal M}.\eqno(3.34)$$
By Lemma 3.3, Eqs.(3.33)-(3.34) imply that
$(Z^2_1-I)PA(I-P)\in{\mathcal Z}({\mathcal M})$. Hence
$(Z^2_1-I)PA(I-P)=0$ holds for all $A\in{\mathcal M}$. Now by Lemma
3.4, one gets $Z^2_1=I$.

{\bf Step 7.}  $\Phi(A)=ZA$ for all $A\in{\mathcal M}$,  where
$Z\in{\mathcal Z}_S({\mathcal M})$ with $Z^2=I$. Therefore, Theorem
3.1 is true.

By Steps 4-6, we have proved that $\Phi(A_{ij})=Z_1 A_{ij}$  for all
$A_{ij}\in{\mathcal M}_{ij}$, where $Z_1\in{\mathcal Z}_S({\mathcal
M})$ with $Z_1^2=I$ ($i,j=1,2$). Define $Z=Z_1$. Now, by a similar
argument to that in the proof of Lemma 2.9, one can show that there
exists a map $f:{\mathcal M}\rightarrow {\mathcal Z}_S({\mathcal
M})$
 such
that $\Phi(A)=ZA+f(A)$ for all $A\in{\mathcal M}$.

To complete the proof of the theorem, we have to prove $f\equiv0$.
Indeed, since $[\Phi(A),\Phi(B)]_*=[A,B]_*$ for every
$A,B\in{\mathcal M}$, we have $[ZA+f(A),ZB+f(B)]_*=[A,B]_*$. It
follows that $Zf(B)(A-A^*)=0$. As $Z^2=I$, we get
$$f(B)(A-A^*)=0\quad{\rm for \ \ all}\quad A,B\in{\mathcal M}.$$
Take $A=iI$ in the above equation, one gets $2if(B)=0$, and so
$f(B)=0$ for every $B\in{\mathcal M}$, completing the proof of
Theorem 3.1. \hfill$\Box$


\end{document}